\documentclass[12pt,a4paper]{article}
\usepackage{graphicx}
\usepackage{amsmath}
\usepackage{amssymb}
\usepackage{amsthm}
\usepackage{geometry}
\usepackage{fancyhdr}%Heads
\usepackage{color} %colors
\usepackage{overpic}

\usepackage[all]{xy}%commutative graph
\newcommand{\be}{\beta}
\newcommand{\ga}{\gamma}

\newcommand{\ra}{\rightarrow}%arrows

\newcommand{\dash}{\textrm{-}}%en dash

\newcommand{\bpf}{\begin{proof}}%thm,prop,lem,pf,example,remark
\newcommand{\epf}{\end{proof}}
\newcommand{\bthm}{\begin{thm}}
\newcommand{\ethm}{\end{thm}}
\newcommand{\bprop}{\begin{prop}}
\newcommand{\eprop}{\end{prop}}
\newcommand{\bcor}{\begin{cor}}
\newcommand{\ecor}{\end{cor}}
\newcommand{\blem}{\begin{lem}}
\newcommand{\elem}{\end{lem}}
\newcommand{\bdefn}{\begin{defn}}
\newcommand{\edefn}{\end{defn}}
\newcommand{\bexmp}{\begin{exmp}}
\newcommand{\eexmp}{\end{exmp}}
\newcommand{\brem}{\begin{rem}}
\newcommand{\erem}{\end{rem}}

\newcommand{\bdia}{\begin{displaymath}\xymatrix}
\newcommand{\edia}{\end{displaymath}}
\newcommand{\beq}{\begin{equation*}\begin{aligned}}
\newcommand{\eeq}{\end{aligned}\end{equation*}}

\newcommand{\intg}{\mathbb{Z}}%numbers
\newcommand{\real}{\mathbb{R}}

\newtheorem{thm}{\textbf {Theorem}}[section]
\newtheorem{cor}[thm]{\textbf{Corollary}}
\newtheorem{prop}[thm]{\textbf{Proposition}}
\newtheorem{lem}[thm]{\textbf{Lemma}}

\theoremstyle{definition}
\newtheorem{defn}[thm]{\textbf{Definition}}

\newtheorem{exmp}[thm]{Example}

\theoremstyle{remark}
\newtheorem{rem}[thm]{Remark}

\newtheorem*{theorem}{\textbf{Theorem}}

\title{Width of the Whitehead double of a nontrivial knot}

\author{Qilong Guo, Zhenkun Li}

\begin{document}
\bibliographystyle{plain}%for reference

\maketitle

%\tableofcontents%table of contents
%\newpage

%————Start from here————

\begin{abstract}    % type your abstract below
In this paper, we prove that $w(K) =4w(J)$, where $w(.)$ is the
width of a knot and $K$ is the Whitehead double of a nontrivial knot $J$.
\end{abstract}

\section{Introduction}
Width is a knot invariant introduced in Gabai \cite{gabai1987foliations}. It has been studied intensively since then. Zupan in \cite{zupan2010properties} conjectured that
$$w(K)\geq n^2 w(J),$$
where $K$ is a satellite knot with companion $J$ and wrapping number $n$. This conjecture is still open. The authors of this paper proved in \cite{guo2018width} a weaker version of the conjecture with wrapping number being replaced by winding number. In this paper we prove a special case of $K$ having wrapping number 2 and winding number 0, which is not covered by the discussion in \cite{guo2018width}.

\begin{theorem}
Suppose $J$ is a nontrivial knot and $K$ is a Whitehead double of $J$. Then we have
$$w(K)=4\cdot w(J).$$	
\end{theorem}

The gap between wrapping number and winding number in the proof in author's previous paper \cite{guo2018width} lies in lemma 4.4, which says that any properly embedded surface representing a generator of $H_2(V,\partial{V})$ has the number of intersection points with the satellite knot $K$ being no smaller than the absolute value of the winding number of $K$. Here $V$ is the tubular neighborhood of $J$ containing $K$. This is in general not true for wrapping number but for the special case that $K$ is a Whitehead double, we can somehow overcome this difficulty. 

{\it\bf Acknowledgement.} 
The first author was supported by NSFC (No.11601519) and Science Foundation of China University of Petroleum, Beijing (No.2462015YJRC034 and No.2462015YQ0604). The second author was supported by his advisor Tom Mrowka’s NSF grant funding No. 1808794.

\section{Width of a Whitehead double}
In what follows, we will use capital letters $K,J,L$ to denote the knot or link classes while use lower case letters $k,j,l$ to denote particular knots or links within the corresponding classes.

\begin{figure}[h]
\centering
\begin{overpic}[width=3.0in]{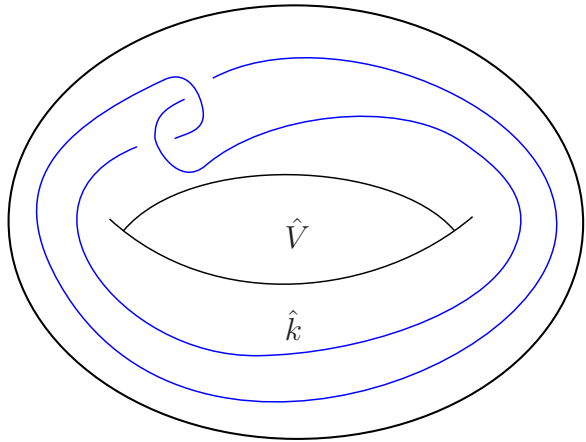}
	\put(48,18){$\hat{k}$}
	\put(48,34){$\hat{V}$}
\end{overpic}
\caption{The pattern of a Whitehead double.}\label{fig_whitehead_double}
\end{figure}

\bdefn\label{defn_whitehead_double}
Suppose $l_w=\hat{k}\cup\hat{j}$ is a Whitehead link in $S^3$. Let $\hat{V}=S^3\backslash{\rm int}({N}(\hat{j}))$ be the exterior of $\hat{j}$ containing $\hat{k}$ in its interior. Since $\hat{j}$ is the unknot in $S^3$, $\hat{V}$ is a solid torus and $\hat{k}\in{\rm int}(\hat{V})$ can be thought  as in figure \ref{fig_whitehead_double}.
Let $j \subset S^3$ be a non-trivial knot and let $V=N(j)$ be the closure of a tubular neighborhood of $j$ in $S^3$. Let $f:\hat{V} \rightarrow S^3$ be an embedding such that $f(\hat{V})=V$, and let $k=f(\hat{k})$. Then $k$ is called a {\it Whitehead double} of $j$ and $j$ is called its {\it companion}.
\edefn

\blem\label{lem_wrapping_number}
Suppose $k,j,V$ are defined as in definition \ref{defn_whitehead_double}, then any meridian disk $D$ of $V$ would intersect $k$ at least two times.
\elem
\bpf
See Rolfsen \cite{rolfsen1976knots}.  
\epf

\brem
The above lemma actually means that the wrapping number of $k$ is 2. On the other hand, $[k]=0\in H_1(V)$ so the winding number of $k$ is 0.
\erem

\bdefn\label{defn_width}
Suppose $h':\real^4\ra\real$ is the projection onto the forth coordinate, $S^3\subset \real^4$ is the unit $3$-sphere, and $h=h'|_{S^3}$. Fix any knot class $K$, let $\mathcal{K}$ be the set of knots $k:S^1\hookrightarrow S^3$ within the class $K$, so that $k$ does not contain the critical points of $h=h'|_{S^3}$, the function $h\circ k: S^1\ra \real$ is Morse, and any two critical points have different critical values.

For any $k\in\mathcal{K}$, let $c_1<...<c_n$ be all the critical values of $h\circ k.$ Pick regular values $r_1,...,r_{n-1}$ so that
$$c_1<r_1<c_2<r_2<...<c_{n-1}<r_{n-1}<c_n,$$
and define
$$w(k)=\sum_{i=1}^{n-1}|h^{-1}(r_i)\cap k|,~w(K)=\min_{k\in\mathcal{K}}w(k).$$
Here $|\cdot|$ means the number of intersection points. Then $w(K)$ is called the {\it width} of $K$.
\edefn

As discussed in the introduction, the key lemma is the following.
\blem\label{lem_planar_surface_intersection}(Key lemma.)
Suppose $S$ is a connected properly embedded planar surface inside $V$ which represents a generator of $H_2(V,\partial{V})$, then $S$ intersects $k$ at least 2 times. 
\elem

We will prove this lemma in next section while first use this lemma to prove our main result.

\bthm
Suppose $J$ is a nontrivial knot (class) and $K$ is a Whitehead double of $J$, then we have
$$w(K)=4w(J).$$
\ethm
\bpf
It is easy to construct examples to show that
$$w(K)\leq 4 w(J).$$
In order to prove the reversed inequality, we repeat the whole argument as in the authors' previous paper \cite{guo2018width} except for applying the key lemma \ref{lem_planar_surface_intersection}. We only sketch the proof as follows. For details, readers are referred to that paper.

First pick a knot $k\in\mathcal{K}$ so that $w(k)=w(K)$. There are corresponding companion knot $j$ and its tubular neighborhood $V$ so that $k\subset V$. Let $T=\partial{V}$. We can perturb $V$ (while fixing $k,j$) so that $V$ does not contain the critical points of $h=h'|_{S^3}$, $h|_{T}$ is Morse and all critical points of $h|_{T}$ have different critical values. 

As in \cite{guo2018width}, there is a Reeb graph $\Gamma_R(V)$ associated to the pair $(V,h|_{V})$. There is a natural way to embed $\Gamma_{R}(V)$ into $V$ and there is a unique loop $l$ in $\Gamma_{V}$ so that

(1). We have that $l$ represents a generator of $H_1(V)$.

(2). We have that $w(L)\geq w(J)$.

(3). We can isotope $l$ into a new position $l'$ with the following significance. If $r$ is a regular level of $h\circ k$, and $m=|h^{-1}(r)\cap k|$, then there exists $P_1,P_2,...,P_m$ as pair-wise disjoint connected components of $h^{-1}(r)\cap V$, so that for each $j$ with $1\leq j\leq m$, $P_j$ intersects $l'$ transversely at one point.

Now let $r_1,...,r_{n-1}$ be all regular values picked as in definition \ref{defn_width}. For each $i$ with $1\leq i\leq n-1$, we pick connected components $P_{i,1},...,P_{i,m_i}$ of $h^{-1}(r_i)\cap V$ as in (3). For each j with $1\leq j\leq m_i$, we know that $P_{i,j}$ is contained in $h^{-1}(r_i)\cong S^2$ so it is a planar surface. From (1) and the fact that $P_{i,j}$ intersects $l'$ transversely at one point, we know that $P_{i,j}$ represents a generator of $H_2(V,\partial{V})$. Hence lemma \ref{lem_planar_surface_intersection} applies and
$$|P_{i,j}\cap k|\geq 2.$$

Now we can apply the identical combinatorial arguments as in the proof of lemma 4.6 in \cite{guo2018width} to finish the proof here.
\epf

%check
\section{The proof of the key lemma}
In order to prove the key lemma, we need some preparations first.

\bdefn\label{defn_auxiliary_function}
Suppose $M=B^3$ or $S^3$ is either a $3$-ball or a $3$-sphere. Suppose $S$ is a properly embedded surface in $M$. This means when $M=B^3$, $\partial{S}\subset \partial{M}$ and when $M=S^3$, $S$ is closed. Then we can define a map
$$C_{M,S}:(M\backslash S)\times (M\backslash S)\ra \{\pm 1\}$$
as
$$C_{M,S}(x,y)=(-1)^{|\gamma\cap S|},$$
where $x,y\in M\backslash S$ are two points and $\gamma$ is an arc connecting two points $x$ and $y$ that is transverse to $S$. The notation $|\cdot|$ means the number of intersection points.
\edefn

The function $C_{M,S}(x,y)$ is independent of the choice of the arc $\ga$. Suppose $\ga$ and $\ga'$ are two such arcs, then $\be=\ga\cup\ga'$ is a closed curve. We know that
$$|\be\cap S|\equiv 0~{\rm mod}~2$$
since $0=[S,\partial S]\in H_2(M,\partial{M})=0$. Hence we conclude
$$|\ga\cap S|\equiv |\ga'\cap S|~{\rm mod}~2.$$

From the definition, we know that $C_{M,S}(x,y)=-1$ would imply that $x$ and $y$ are not in the same component of $M\backslash S$. 

Another good property of this function is the following equality: for any 3 points $x,y,z\in M\backslash S$, we have
\begin{equation}\label{eq_composition}
	C_{M,S}(x,y)\cdot C_{M,S}(y,z)=C_{M,S}(x,z).
\end{equation}

\brem
This type of functions are called potentials in Gabai \cite{gabai1987foliations}.
\erem

\blem\label{lem_nonseparating}
Suppose $S$ is a connected, properly embedded planar surface inside a 3-ball $B$. Let $B_1$ and $B_2$ be the two components of $B\backslash S$. Then there is no Hopf link $l_h=l_1\cup l_2$ inside $B$ such that 
$$l_i\subset B_i~for~i=1,2.$$
\elem
\bpf
Assume, on the contrary, that there is a Hopf link $l_h=l_1\cup l_2$ such that $l_i\subset B_i.$

Note $B\backslash S$ has two components. To show this, let $N(S)$ be a product neighborhood of $S$, and $B'=B\backslash {\rm int}(N(S))$. We know that ${\rm int}(B')$ and $B\backslash S$ are homotopic. From the connectedness of $S$ and Mayer-Vietoris sequence
\begin{equation*}
\xymatrix{
H_1(B)\ar[r]\ar@{=}[d]&H_0(\partial{N(S)})\ar[r]\ar@{=}[d]&H_0(N(S))\oplus H_0(B')\ar@{=}[d]\ar[r]&H_0(B)\ar@{=}[d]\ar[r]&0\\
0\ar[r]&\intg\oplus\intg \ar[r]&\intg\oplus H_0(B')\ar[r]&\intg \ar[r]&0
}	
\end{equation*}
we know that $H_0(B')=\intg\oplus\intg$. Let $x_1\in l_1$ and $x_2\in l_2$ be two points, then from the hypothesis of the lemma, $x_1$ and $x_2$ lie in two different components of $B\backslash S$. Since there are only two components of $B\backslash S$, we know that
$$C_{B,S}(x_1,x_2)=-1,$$
and hence there exists a curve $\ga$ connecting $x_1$ with $x_2$ and having an odd number of intersections with $S$.

We can pick another $3$-ball $B''$ and glue it to $B$ along their common spherical boundary to get an $S^3$. Inside this $S^3$, $l_h=l_1\cup l_2$ is still a Hopf link. Inside $B''$, we can pick disks $D_1,...,D_t$ so that 
$$S^2=S\cup D_1\cup,...,D_t\subset S^3$$
is a $2$-sphere. Since $\ga\subset B$, we know that $\ga\cap S^2=\ga\cap S$ and conclude that
$$C_{S^3,S^2}(x_1,x_2)=-1,$$
so the two components of the Hopf link $l_h$ lie in two different components of $S^3\backslash S^2$, which is absurd since it is well-known that the Hopf link is non-separating.
\epf

Now we are ready to prove the key lemma.

\bpf[proof of lemma \ref{lem_planar_surface_intersection}.]
Suppose $k$ is a Whitehead double of $j$, $V$ is the tubular neighborhood of $j$ containing $k$ and $S$ is a connected, properly embedded planar surface in $V$, representing a generator of $H_2(V,\partial V)$. Assume that $S$ and $k$ have less than 2 intersections, then since $[k]=0\in H_1(V)$, we know that $S$ and $k$ must be disjoint.

\begin{figure}[h]
\centering
\begin{overpic}[width=4.0in]{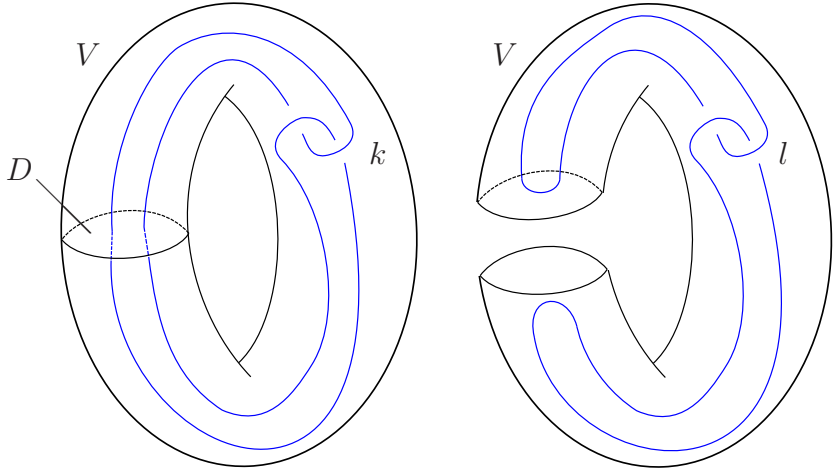}
	\put(-3,38){\line(1,-1){7}}
	\put(-7,38){$D$}
	\put(40,40){$k$}
	\put(93,40){$l$}
	\put(56,53){$V$}
	\put(2,53){$V$}
\end{overpic}
\caption{After cutting and re-guling we will get a Hopf link $l$.}\label{fig_hopf_link}
\end{figure}

It is straightforward to see that there is a meridian disk $D$ of $V$, such that in $B=V\backslash {\rm int}(N(D))$, $D$ intersects $S$ transversely, and that after adding two small arcs to $(k-{\rm int}({N}(D)))$ near $\partial{N(D)}$, we will get a Hopf link $l$ out of $k$. Here $N(D)$ is a neighborhood of $D$ in $V$. See figure \ref{fig_hopf_link}.

\begin{figure}[h]
\centering
\begin{overpic}[width=4.5in]{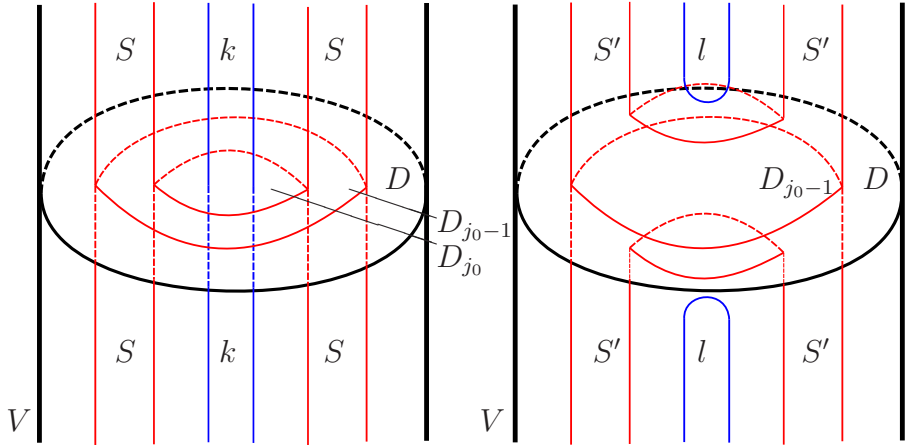}
	\put(40,30){$D$}
	\put(36,30){\line(3,-1){10}}
	\put(27,30){\line(3,-1){19}}
	\put(46,25){$D_{j_0-1}$}
	\put(46,21){$D_{j_0}$}
	\put(95,30){$D$}
	\put(83,30){$D_{j_0-1}$}
	\put(21,45){$k$}
	\put(9,45){$S$}
	\put(33,45){$S$}
	\put(9,10){$S$}
	\put(33,10){$S$}
	\put(76,45){$l$}
	\put(64,45){$S'$}
	\put(88,45){$S'$}
	\put(64,10){$S'$}
	\put(88,10){$S'$}
	\put(21,10){$k$}
	\put(76,10){$l$}
	\put(-3.5,2){$V$}
	\put(51,2){$V$}
\end{overpic}
\caption{Surgery on $S$ along $D_{n}$, ..., $D_{j_0-1}$.}\label{fig_surgery_on_s}
\end{figure}

If $D\cap S=\emptyset$, then we can apply lemma \ref{lem_nonseparating} directly to conclude a contradiction. Note the two components of the Hopf link $l$ are in two different components of $B\backslash S$. This argument for this is exactly the same as the following more general case (when $D\cap S\neq \emptyset$): in figure \ref{fig_surgery_on_s}, we have
$C_{B,S}(u,v)=-1$ for the same reason while $C_{B,S}(x,u)=C_{B,S}(y,v)=0$ as there is no need for doing any surgeries.

If $D\cap S\neq \emptyset$, we can assume that 
$$D\cap S=\be_1\cup \be_2...\cup\be_m,$$
where $\be_1,...,\be_n$ are intersection circles of the two surfaces and are in the order such that $\be_j$ bounds a disk $D_j\subset D$ disjoint from any $\be_i$ for $i<j$.

If all $D_i$ are disjoint from $k$, then we can do a series surgeries on $S$ with respect to $D_n,D_{n-1},...,D_{1}$ one by one to get a surface $S''$ so that $S''$ is disjoint from $D$ and $k$. We can pick some connected component of $S''$, and it will also have such properties that we can apply the argument above to get a contradiction.

Now we are in the most complicated case where some $D_i$ intersects $k$. Suppose $j_0$ is the greatest index such that $D_j\cap k\neq\emptyset$. We first claim that $D_{j_0}$ cannot have a unique intersection point with $k$. Suppose the contrary, then a sequence of surgeries on S with respect to $D_n,...,D_{j_0+1}$ would generate a surface $S'$ such that $S'$ represent a generator of $H_{2}(V,\partial V)$, $S'$ is disjoint from $k$ and
$$D_{j_0}\cap S'=\partial D_{j_0}=\be_{j_0}.$$
Suppose $S'_0$ is the component of $S'$ containing $\be_j$, then $S_0'$ is still a planar surface. A surgery on $S'_0$ with respect to $D_{j_0}$ would result in two surfaces $S_1'$ and $S_2'$, each of which contains one copy of $D_{j_0}$ and hence has a unique intersection point with $k$. This is impossible since $[k]=0\in H_1(V)$.

Hence we conclude that $D_{j_0}$ has two intersection points with $k$, which are all the intersection points of $D$ with $k$. Now, as before, we can do a sequence of surgeries on $S$ with respect to $D_n,...,D_{j_0+1}, D_{j_0}$ to get $S'$. When doing last surgery with respect to $D_j$, we shall modify $k$ at the same time: cut $k$ along $D$ and glue two small arcs to the newly born boundary points near $D$ to get a Hopf link $l$ disjoint from $S'$ and $D$. See figure \ref{fig_surgery_on_s}. Now we can apply the remaining sequence of surgeries on $S'$, with respect to $D_{j_0-1},...,D_1$, to get a surface $S''$ which represents a generator of $H_{2}(V,\partial V)$ and is disjoint from $D$ and $l$.

\begin{figure}[h]
\centering
\begin{overpic}[width=2.5in]{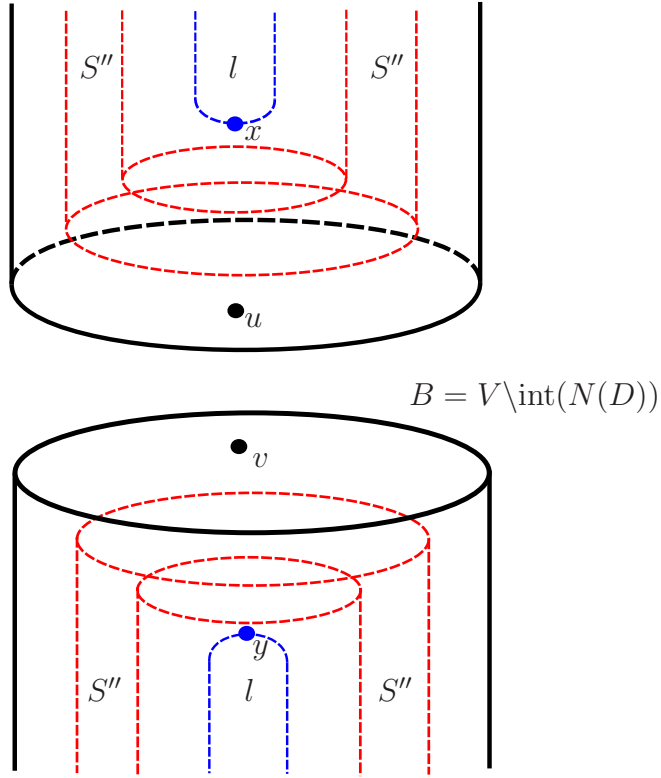}
	\put(30,10){$l$}
	\put(28,90){$l$}
	\put(9,90){$S''$}
	\put(46,90){$S''$}
	\put(10,10){$S''$}
	\put(47,10){$S''$}
	\put(30,82){$x$}
	\put(30,58){$u$}
	\put(31,40){$v$}
	\put(31,16){$y$}
	\put(51,48){$B=V\backslash{\rm int}({N}(D))$}
\end{overpic}
\caption{Surgery on $S$ along $D_{n}$, ..., $D_{j_0-1}$.}\label{fig_four_points}
\end{figure}

Now we can cut off a neighborhood $N(D)$ of $D$ from $V$, which is disjoint from $l$ and $S''$, to get a $3\dash$ball $B=V\backslash {\rm int}({N}(D))$. In order to see that the two components of $l$ are in different components of $B\backslash S''$, we can pick points $x,y$ on two arcs of $l$ which do not belong to k, and pick two points $u,v$ near the boundary $\partial{N(D)}$ but on different sides of $D$. See figure \ref{fig_four_points}.

%\begin{figure}[h]
%	\center
%	\includegraphics[width=2in]{4.jpg}
%	\caption{Points $u,v,x,y$}\label{fig_four_points}
%\end{figure}

There is a symmetry between the two pairs $(x,u)$ and $(y,v)$ with respect to $D$ so we have
$$C_{B,S''}(x,u)=C_{B,S''}(v,y),$$
where $C_{B,S''}$ is defined as in definition \ref{defn_auxiliary_function}. Since $S''$ represents a generator of $H_{2}(V,\partial V)$, we have
$$C_{B,S''}(u,v)=-1.$$
Using equality (\ref{eq_composition}), we have
$$C_{B,S''}(x,y)=C_{B,S''}(x,u)\cdot C_{B,S''}(u,v)\cdot C_{B,S''}(v,y)=-1.$$
Thus the two components of $l$ are in different components of $B\backslash S''$. Pick a connected component $S_0''$ which separates the two components of $l$, we can apply lemma \ref{lem_nonseparating} to get a contradiction.
\epf

%————End from here————

\newpage
\bibliography{Index}%for reference
\end{document}